\documentclass[11pt]{article}

\input amssym.def
\input amssym.tex

\newtheorem{theorem}{{Theorem}}[section]
\newtheorem{isom.ext}[theorem]{{Trivial isometric extension}}
\newtheorem{lemma}[theorem]{{Lemma}}
\newtheorem{corollary}[theorem]{{Corollary}}
\newtheorem{fact}[theorem]{{Fact}}
\newtheorem{remarks}[theorem]{{Remarks}}

\newenvironment{proof}{{\it Proof.}}{\hfill$\diamondsuit$\medskip}

\begin{document}

\title{On Lorentz dynamics : From group actions to warped products via homogeneous spaces}
\author{A. Arouche, M. Deffaf and A. Zeghib}
\date{\today}
\maketitle

\begin{abstract}
We show a geometric rigidity of isometric actions of non compact
(semisimple) Lie groups on Lorentz manifolds. Namely, we show that
the manifold has a warped product structure of a Lorentz
 manifold with constant curvature by a Riemannian manifold.
\end{abstract}

\section{Introduction}
Recall the following result of \cite{ze1}, which shows how homogeneous spaces are
rare in Lorentz geometry (in comparison with the Riemannian case, say)

\begin{theorem} \label{zeghib.symmetric} \cite{ze1}Let $(M, g)$ be  a   homogeneous Lorentz
 space of dimension
$\geq 3$, with {\em  irreducible} isotropy group,
then it has constant sectional curvature.
\end{theorem}

Observe that the statement in \cite{ze1} seems weaker, since the isotropy group is assumed to satisfy
a supplementary condition: non-precompactness.  However, this follows from irreducibility.  Indeed, in the same vein as \cite{ze1}, the principal result of \cite{B-Z} says how irreducibility is strong in
the Lorentz setting:
\begin{theorem}  \label{boubel.zeghib} \cite{B-Z} A  Lie subgroup (not assumed a priori to be closed)  of $O(1, n)$,  which does not preserve any one dimensional
isotropic subspace of ${\Bbb R}^{1+n}$, is up to conjugacy,
a union of some components of some $O(1, p) \subset O(1, n)$.

\end{theorem}

Our goal in the present article is to relax homogeneity by considering (non-transitive) isometric group actions. This work is actually motivated by the study of isometric Lie group actions on {\it non-compact} Lorentz manifolds, for instance in the same vein as \cite{ada, Kow1, Kow2}...

\subsection{Warped product structure versus partial homogeneity}
 We ask  firstly if there is an adaptation
of  Theorem \ref{zeghib.symmetric} to non-transitive isometric actions.
 In this situation, we consider a group $G$ acting isometrically
 on a Lorentz manifold $(M, g)$. Each orbit is a homogeneous space. However, the
causal type of the orbit may be, timelike, spacelike or lightlike, that is, the induced metric
is Lorentzian, Riemannian or degenerate, respectively.
The following  generalization  of Theorem \ref{zeghib.symmetric}
relies on the existence of orbits of Lorentz type satisfying
irreducibility.  It says roughly  that the space is partially of
constant curvature.

\begin{theorem} \label{non.transitive}
Let  $G$ be a Lie group acting isometrically on  a   Lorentz manifold  $(M, g)$
of dimension $\geq 3$. Suppose there exists an orbit $N$ which is
a (homogeneous) Lorentz space
with  irreducible isotropy.

Then, $N$ has constant (sectional) curvature, and a neighborhood
of it  is a warped product $L \times_w N$,   where $L$ is some
Riemannian manifold. Furthermore, the factor $N$ corresponds to
the orbits of $ G$.

 \end{theorem}

Definition and fundamental properties  of warped products are in \S \ref{section.non.transitive}.
It follows that the group $G$ is a subgroup of the isometry group of a
 constant curvature manifold
$N$ (acting transitively on it). It is a non-difficult algebraic matter to classify them. Conversely, any such group acts isometrically on any warped product
 $L\times_w N$.


\subsection{Non-properness  versus Irreducibility}
 Let us go a step further, and try to get rid of the irreducibility hypothesis.  In fact, irreducibility
  is an algebraic condition which looks somehow non-adapted to our dynamico-geometrical setting here.  We  want to
 substitute for it a more natural  dynamical condition. Our theory is that non-properness is
 good enough for this role.

 \subsubsection{Recalls}
  We find it worthwhile to make some order around the concept of non-properness of actions.
   This will be useful in the sequel (statements and proofs).

Recall that an action of a group $G$ on a space $M$ is called {\bf proper}, if for any
sequences $(x_n)$ of $M$, and $(g_n)$ of $G$, if $(x_n)$ and $(g_nx_n)$ converge in $M$, then
some subsequence of $(g_n)$ converges in $G$.

For our purpose here the following variant will be useful. We say
that the action of $G$ is {\bf locally equicontinuous}, if keeping
the notations as above, a subsequence of $(g_n)$ is (locally)
equicontinuous (or say it is equicontinuous in a neighborhood of
the limit of $(x_n)$). Therefore,  a subsequence of $(g_n)$ is
converging in the group of homeomorphisms of $M$, but, the limit
does not necessarily belong to $G$.
 Obviously, a non locally-equicontinuous action is non-proper. The converse
 is not true. The standard example of a non-proper but equicontinuous action is the usual linear
action of ${\Bbb R}$ on the torus with dense orbits. This is in general the case of any
non-closed Lie group of the isometry group of a compact Riemannian manifold.
Another example is the action of the universal cover $\tilde{G}$ on $G$
(via the canonical projection). It is always locally equicontinuous, but proper only if
$G$ has a finite fundamental group.
Observe nevertheless:
\begin{fact} Let $G$ be a Lie group acting by preserving a pseudo-Riemannian structure on a manifold $M$. If $G$ is the full isometry group, or $G$ is
 semi-simple with finite center, then its action is non-proper iff it is non-equicontinuous.
\end{fact}
For the proof, recall the well  known fact that a $C0$-limit of
pseudo-Riemannian (smooth) isometries is a smooth isometry, and
that the Lie group topology coincides with the $C0$ topology.
This is equivalent to saying  that the isometry group is closed in
the group of homeomorphisms. (Actually, this fact is  general for
all rigid geometric structures).
For $G$ a simple Lie group with finite center, recall that its
image under a homomorphism into any Lie group is closed, and that
$G$ is a finite cover of it. An analogous argument applies to the
semi-simple case with finite center. \hfill $\diamondsuit$

Therefore, in statements (essentially inside proofs) below,which
involve  semi-simple Lie groups, we will not worry to sway from
compactness to pre-compactness.

 A $G$-homogeneous space $G/H$ is non-proper if the $G$-left action
 on it is. This is equivalent, in the general case,  to the fact  that $H$ is not precompact, and
 to that $H$ is not compact in the semi-simple case.

 \subsubsection{Semi-simple group actions with non-proper orbits}
 Without   a priori irreducibility hypothesis, we  have
 the following generalization of Theorem \ref{non.transitive}, assuming
 the orbits are non-proper, and
   the group $G$ is semi-simple (a kind of intrinsic
 irreducibility).

 \begin{theorem} \label{semisimple}
Let  $G$ be a    {\em semi-simple}Lie group of finite center
acting isometrically on  a   Lorentz manifold  $(M, g)$ of
dimension $\geq 3$.
Suppose  that no  (local) factor of $G$ is locally isomorphic to $SL(2, {\Bbb R})$
 and that there exists a non-proper  orbit $N$ of  Lorentz type,
(that is   $N$ has   a  non-compact  isotropy).

 Then, up to a finite cover, $G$ factorizes $G = G_1 \times G_2$, where:

 - $G_2$ acts properly on a neighborhood of $N$, with spacelike orbits of constant dimension.

 -  On a $G$-invariant neighborhood  $U$ of $N$,  $G_1$ acts isometrically on the Lorentz quotient $U/G_2$
which is a warped product
$L \times_w N_1$,
where $N_1$ has constant (sectional) curvature
 and corresponds to the $G_1$-orbits. In particular, if $G$ is simple, then $U$ itself
 is a warped product.

\end{theorem}
Looking at semi-simple Lie groups acting transitively, with
non-compact isotropy on a constant curvature Lorentz space, on can
prove in a  standard way the following:

\begin{corollary} In the case above, up to a (central) cover, $G_1$ is
$O(1, n)$ (resp.  $O(2, n)$)   the isometry group of the de Sitter space (resp. anti de Sitter space), that is the universal Lorentz space
of positive (resp. negative) curvature.

\end{corollary}

\begin{remarks} {\em   ${}$

1) In both theorems above, the warped product is local, i.e. not
the whole space is a warped product.  To see this, one considers
the $O(1, n)$-action on the Minkowski space ${\Bbb R}^{1, n}$. If
the Lorentz  quadratic form is $q= -x_02 + x_12 + \ldots +
x_n2$, then the warped product is defined exactly on the region
${ q >0}$.

2) The result does not seem to be optimal, that is, it might be generalized to
other groups.

3) Warped product structures on universal covers of {\it compact} Lorentz manifolds with strong dynamics, are
obtained, for instance, in \cite{Gro, ze2, ze3}.
}
\end{remarks}

\subsection{From non-proper actions to non-proper homogeneous spaces }

 Let us go another step,   by asking how to get such non-proper orbits from a
 global condition on the action? For instance, is  orbital non-properness
 inherited from non-properness of the (ambient) action?

\begin{theorem} \label{action.homogeneous}
Let $G$ be a semi-simple Lie group  of finite center acting
isometrically and {\em nonproperly}  on a Lorentz manifold $M$.
Suppose  that no  (local) factor of $G$ is locally isomorphic to
$SL(2, {\Bbb R})$. Then,  there
 is a point  with a non-compact stabilizer.   In particular,
 the restriction of the action of $G$ to its orbit   is nonproper.
  More exactly, the stabilizer of some point contains a non-trivial unipotent one-parameter group.
(In other words,  a non proper Lorentz $G$-space contains  a non-proper
$G$-homogeneous orbit)
\end{theorem}

This result allows one to get from (non-transitive) actions to
homogeneous (i.e. transitive) ones. This is a  common philosophy
for actions with strong dynamics and a geometric background.   The
result here is  in particular reminiscent to
  the so-called  Zimmer's embedding Theorem (see for instance
\cite{Zim}). ``Unfortunately'', there is a damper to put on : the
orbit is a non-proper homogeneous space, but not necessarily
Lorentz! The nuisance
  is that it can be lightlike (degenerate);  another story.

\section{Proof of Theorem \ref{non.transitive}} \label{section.non.transitive}

\subsection{An algebraic lemma}
\begin{lemma}  \label{algebraic}
 Let $E$ (resp. $F$)  be a Lorentz (resp. Euclidean) vector space. Denote by
$O(E)$ and  $O(F)$    their respective orthogonal groups.
Let $H^\prime$ be a Lie subgroup of $O(E)  \times O(F)$, whose projection
on $O(E)$   acts irreducibly on $E$.

 Then,  $H^\prime$
contains a subgroup $H \subset O(E) \times \{ 1\}$, which contains
the identity  component of $O(E)$. In particular:

- Any  linear $H$-invariant  mapping $f: E \to F$ ($f \circ h = f$, for
any $h \in H$) is trivial.

- The same is true for any $H$-invariant  bilinear antisymmetric
mapping $E \times E \to F$.

\end{lemma}

\begin{proof}
The crucial fact  follows from Theorem  \ref{boubel.zeghib}
applied it  to $H$,  the projection of $H^\prime$ on $O(E)$.  It
is a finite  union of components  of $O(1, p)$. Say $H = O(1, p)$
to simplify notation.    Since $O(F)$  is compact, $H$ is
isomorphic to  the non-compact semi-simple Levi factor of
$H^\prime$. Therefore, (up to a cover...) $H^\prime$ contains a
subgroup isomorphic to $H$, that is, there exists a homomorphism
$\rho: H= O(1, p) \to O(F)$, such that the graph $\{ (h, \rho(h)),
h \in O(1, p) \}$ is contained in $H^\prime$.

Next, one checks $\rho$ is trivial. This    uses a basic fact of Lie groups theory: a semi-simple Lie group
of non-compact type has no non-trivial homomorphism into a compact group.  The idea, in this case, is that in $O(1, p)$ there are one parameter groups (the unipotent ones)  having all their non-trivial
elements conjugate (this is easy to see in the case of $O(1, 2)$ which is essentially
$PSL(2, {\Bbb R})$).  Such a conjugacy is impossible in a compact group.

For the last two conclusions of the lemma, one can assume $F =
{\Bbb R}$. The kernel of th linear mapping $f$ is $O(1,
p)$-invariant; hence it is trivial by irreducibility. A similar
argument yields triviality of invariant antisymmetric bilinear
mappings.
\end{proof}

\subsection{Group actions}

\begin{lemma} \label{integrability} Let $G$ be a Lie group acting isometrically on a Lorentz manifold
$(M, g)$. Let $N$ be an orbit of $G$, which is of Lorentz type and has
an irreducible isotropy group (inside $G$).

Then, the same is true for all orbits in a neighborhood of $N$. In
particular, the orbits of $G$ determine a foliation (i.e. have a
constant dimension). Furthermore, the orthogonal distribution of
this foliation is integrable.
\end{lemma}

\begin{proof} Consider  $x_0 \in N$, and denote by $H$ its isotropy group. The orthogonal space
$L_{x_0}$ of $T_{x_{0}}N$ in $T_{x_{0}}M$  is spacelike (the
metric on it is definite positive). We are in position to apply
Lemma \ref{algebraic} with $E= T_{x_0}N$ and $F= L_{x_0}$.
It then follows that the action of $H$ on $L_{x_0}$ is trivial. Let $\exp_{x_0}$ denote
  the exponential of the Lorentz metric and consider the (local) submanifold
${\cal L}_{x_0} = \exp_{x_0} (L_{x_0})$. Then $\exp_{x_0}$
conjugates the infinitesimal action of $H$ on $L_{x_0}$, with its
action on ${\cal L}_{x_0}$. In particular, $H$ acts trivially on
this latter submanifold.  That is $H$ is contained in the isotropy
group of any point of ${\cal L}_{x_0}$.  An obvious
semi-continuity argument implies that isotropy groups can not be
bigger.
Therefore, we have a foliation by $G$-orbits, all satisfying  the
same irreducibility condition for their isotropy groups. Let us
denote this foliation by ${\cal N}$ and its tangent bundle by
$T{\cal N}$.
Let $L$ be the orthogonal distribution.  The obstruction to integrability of $L$
can measured by means of a   tensor $T: L \times L \to T{\cal N}$. It is
defined by $T(X, Y)= $  the orthogonal projection on $T{\cal N}$ of the bracket
$[X, Y]$, where $X$ and $Y$ are sections of $L$.  Since the isotropy group
acts trivially on $L$ and irreducibly on $T{\cal N}$,  $T$ is trivial, that is
$L$ is integrable.
\end{proof}

\subsection{Warped product, end of the proof}
\label{warped}

Let   $(L, h)$ and
$(N, m)$ be   two pseudo-Riemannian
  manifolds
 and $ w: L \to  {\Bbb R}^+-\{0\}$
a {\it warping} function.
  The warped
product \mbox{$M= L\times_wN$}, is the topological
product
  $L \times N$, endowed with
the  pseudo-Riemannian  metric
$g= h \bigoplus w m$.

Our goal now is to prove that $M$ is a warped product. So far, we
have the orthogonal foliations ${\cal N}$ and ${\cal L}$.
 One can say that
De Rham decomposition theorem is a criterion for a couple of such
foliations in order that they determine a (local) {\it direct}
pseudo-Riemannian product. The condition is that (the tangent
bundles of) ${\cal N}$ and ${\cal L}$ are parallel, or a priori
more weakly, that leaves of ${\cal N}$ and ${\cal L}$ are
geodesic. There is a  similar, but more complicated,  criteria for
warped products \cite{Hie, P-R}. We will not use this criterion,
but rather give a brief proof in our case. Our terminology here is
close to that of \cite{ze3}, which may be consulted for a more
complete exposition.
Let $N$ and $L$ be (local) leaves of a point $x_0$ for the
foliations ${\cal N}$ and ${\cal L}$ respectively.  So, locally,
$M$ has an adapted topological product $L \times N$. The metric
can be written
 $$g_{(l, n)} = h_{(l, n)} \bigoplus  m_{(l, n)}$$

$\bullet$ Let us show that $h_{(l, n)} = h$, that is,  it does not
depend on $n$. This is clear since $G$ acts isometrically: if $k
\in G$, then it  sends ${\cal L}_{(l, n)}$ to ${\cal L}_{k(l,
n)}$, where $k(l, n)$ has the form $(l, n^\prime)$ (orbits of $G$
correspond to ${\cal N}$). Therefore $g=  h \bigoplus m_{(l, n)}$
(the geometric meaning of this fact is   that ${\cal N}$ is a
geodesic foliation \cite{ze3}).

 $\bullet$ In order to understand the variation of $m_{(l, n)}$ as a function
 of $(l, n)$, write  $x_0 = (l_0, n_0)$, fix  $l_1 \in L$  and consider the mapping
 $$S: (l_0, n) \in N = {\cal N}_{(l_0, n_0)} \to  (l_1, n) \in {\cal N}_{(l_1, n_0)}$$
 $S$ commutes with the $G$ action on the $G$-orbits of $(l_0, n_0)$ and $(l_1, n_0)$.
 In particular it commutes with the isotropy actions at these two points. As showed previously these isotropy groups are the full orthogonal groups
of the Lorentz scalar products on their tangent spaces. In
particular, they preserve, up to a multiplicative constant, only
one Lorentz scalar product. This means that $S$ is a homothety at
$(l_0, n_0)$ : the image metric  equals the metric at $(l_1, n_0)$
(along $ {\cal N}_{(l_1, n_0)}$) up to a multiplicative factor
$w(l_0, n_0)$. Now, since $S$ commutes with the (full) action on
orbits, it follows that $w$ does not depend of $n$. That is, if
$m= m_{(l_0, n)}$ is the metric on $N$, then $m_{(l, n)} = w(l)
m$. In sum, $g= h \bigoplus w(l) m$, that is,  $M$ is a warped
product.

$ \bullet$ Finally,  it remains to see that $N$ has constant curvature.  This is exactly the content of Theorem \ref{zeghib.symmetric} since $N$ has a non-precompact irreducible
isotropy.
\hfill  $\diamondsuit$

\section{Proof of Theorem \ref{semisimple}}

We will in fact prove Theorem \ref{semisimple} under the homogeneity assumption,  that
is $G$ acts transitively on $M$. This will be
a generalization Theorem \ref{zeghib.symmetric}, where one keeps
  non-precompactness assumption together with semi-simplicity of the group, and gives up the irreducibility one.
  The proof in the non transitive case will be  just a variation of that of Theorem \ref{non.transitive} using the transitive statement which is:

 \begin{theorem} \label{semisimple.bis}
Let $(M, g)$ be  an {\em irreducible}   $G$-homogeneous Lorentz
 space of dimension
$\geq 3$, with  {\em non-precompact} isotropy group,  and
$G$  a semi-simple Lie group with
no  (local) factor  locally isomorphic to $SL(2, {\Bbb R})$.

$\bullet$ Then, the isotropy group is irreducible and
$M$   has constant sectional curvature.

- In the general case where $M$ is not assumed to be irreducible, we have:

$\bullet$ $M$ is locally a direct product $M = L \times N$ where $L$ is  a Riemannian
homogeneous manifold, and $N$ is Lorentz and has constant curvature. To this
splitting corresponds an analogous (local) one for $G$.

\end{theorem}

\begin{proof}
   For $x \in M$, and a lightlike (i.e. isotropic) vector $u \in T_xM$,
consider the orthogonal  hyperplane $ u^\perp$.  Let $C_x$ the set
of those $u$, for which $u^\perp$ is tangent to a totally geodesic
(lightlike) hypersurface, that is,
$\exp_x (u^\perp)$ is a totally geodesic hypersurface (near $x$). The crucial fact, proved
in \cite{ze1} is that non-precompactness of the isotropy group $H_x$
implies $C_x$ is non-empty.

-  Suppose $C_x$ is finite.  One can   (locally) define only finitely
many  continuous  sections
$x \to u(x) \in C_x$.  In particular, one can suppose these sections invariant under
the $G$-action. In fact, to simplify notation in the following argument, one is allowed
to suppose that $C_x$ has (everywhere) cardinality 1. Therefore, we have a $G$-invariant
distribution of hyperplanes $x \to u(x)^\perp$. It is integrable, the leaf at $x$ being  the geodesic hypersurface $$ {\cal H}_u= \exp_x (u^\perp)$$
We get from this that $M$ possesses
a codimension one $G$-invariant foliation. The quotient space is a 1-manifold.  But
a simple Lie group acting  (non-trivially) on a 1-manifold must be
 locally isomorphic to $SL(2, {\Bbb R})$. This is impossible because of our assumption
 on $G$.

--  It then follows that $C_x$ is infinite.  One then shows in a standard way that
there must exist  a subspace $E_x$ which is:  generated by $E_x \cap C_x$,
spacelike,  and on which  the isotropy group is irreducible.

$\bullet$  Therefore, we have a distribution $E$ on which the
isotropy group acts irreducibly.   As in the proof of Theorem
\ref{non.transitive}, one defines an integrability obstruction
tensor for $E$, which must vanish, by Lemma \ref{algebraic}  since
it is antisymmetric and invariant under the isotropy group.
Therefore $E$ is integrable. We denote by ${\cal N}$ its tangent
foliation.

$\bullet$ Also   $L = E^\perp$ is integrable. Indeed, the leaf of $L$ at $x$ is nothing but
the intersection $${\cal L}_x = \cap_{u \in C_x} {\cal H}_u$$
In addition, as all the hypersurfaces ${\cal H}_u$ are geodesic, the foliation ${\cal L}$
is geodesic.  We recalled  in \S \ref{warped} the interpretation of being geodesic by the fact that the holonomy mappings of the orthogonal foliation ${\cal N}$, defined as mappings between (local)
leaves
  ${\cal  L}$, are {\it isometric}.

$\bullet$ A leaf $ N$ of ${\cal N}$ is a Lorentz  manifold with a
big, in fact maximal isotropy group (at each point), that is
(essentially) $O(1, p)$ (where $p+1 = $ dim$M$). At this stage, we
don't know if $M$ is homogeneous. However, remember the proof in
\cite{ze1}:  non-precompacteness and irreducibility of isotropy
groups lead to existence of geodesic hypersurfaces, which leads in
turn to that the sectional curvature is constant for all
$2$-planes tangent to a same point. Then, using Schur's lemma one
proves $N$ has (everywhere the same) constant curvature.

$\bullet$ The isotropy group of any point $x \in N$ acts trivially
on the leaf ${\cal L}_x$. It follows that the group $R$ generated
by all the isotropy groups of points of $N$, preserves
individually the leaves of ${\cal N}$. It is known, that for
constant curvature spaces, the isotropy group of two different
points generate the full isometry group. In particular $R$ acts
transitively on the leaves of ${\cal N}$.
In fact, the  foliation ${\cal N}$ is defined by the $R$-action.
 We are thus exactly in the situation of Theorem
  \ref{non.transitive} (where $R$ plays the role of $G$).
   Therefore, we deduce that $M$ is (locally) a warped product $L \times_w N$.

$\bullet$ The quotient $M/ {\cal L}$ has a similarity Lorentz
structure, that is, a Lorentz metric up to a (global) constant,
preserved by  $G$. On  other words, $G$ acts by homothety on the
constant curvature $N$.  It is easy to describe the similarity
group of the Minkowski space. One can in particular see that a non
semi-simple Lie group acts {\it transitively} by homothety.

 $\bullet$ We infer from this that $N$ has a non-vanishing curvature. Since $G$ acts
 transitively on $M = L \times_w N$ by preserving the warped product structure, all
 the leaves $\{l \} \times N$ are isometric, and hence have a same curvature. However,
 metrics at two levels $l_1$ and $l_2$ are related by a factor $ \frac{w(l_1)}{w(l_2)}$. Curvature are related by the inverse ratio.  From constancy
 of curvature, we infer that $w$ is a constant function, that is $M= L \times N$ is a
 {\it direct} product. This finishes the proof of Theorem \ref{semisimple.bis}.
  \end{proof}

\section{Proof of Theorem \ref{action.homogeneous}}

The  following method has become a standard  ingredient
in the study of ``geometric'' $G$-actions, see for instance
\cite{ada, A-S, Kow1, Kow2}...
One considers the action of the group
$G$ on the space  $S^{2}(\cal{G})$ of symmetric bilinear forms on
its  Lie algebra $\cal G$.
There is a  Gauss  $G$-equivariant  map $\Phi: M \to S^{2}(\cal{G})$.
 Non-properness of the $G$-action on $M$ translates
 to a  non-properness of
 the action of $G$ on the image $\Phi (M)$. This latter is ``algebraic'',   it has
 a poor  dynamics,   easy to understand.  From this, one hopes to get information about
 the $G$-action on $M$.

In our case here, one shows there
exists a point $q\in M$ such that $\cal G$ admits an isotropic
subspace with respect to the symmetric bilinear form  $\Phi (q)$,
 of dimension $\geq 2$. Then the non-precompactness of the stabilizer $stab(q)$ follows.

\subsection{The Gauss map}
Let $G$ be a Lie group acting by isometries  on a Lorentz manifold
  $(M, h)$. For each $X\in \cal{G}$ let $\overline {X}$ be the
vector field on $M$ given by:
\[\overline {X}_{x}=\frac{d}{dt}(exp(tX).x)|_{t=0}.\]
Let $\Phi:M\rightarrow S^{2}(\cal G)$ be the so-called Gauss map
given by \[x\mapsto\Phi_{x}:(X,Y)\mapsto h_{x}(\overline
{X_{x}},\overline {Y_{x}}).\]

Recall the definition of the $G$-action on $S^{2}\cal G$.  It is given by :
\[ (g.q)\left(X_{1}, X_{2}\right)=q \left(Ad_{g^{-1}}X_{1},\ Ad_{g^{-1}}X_{2}\right)\].
for $q$ in $ S2\left(\cal G\right)$ and $g$ in $G$, and $X_{1}, X_{2}\in \cal G$.

Then $\Phi$ is equivariant,  that is : \[g.\Phi_{x}=\Phi_{g.x}
\forall g \in G. \] Indeed, for $g \in G$ and $X \in \cal G$, we have :
\[\overline {Ad_{g}X_{gx}}=\frac {d}{dt}(exptAd_gX.gx)\mid _{t=0}=\frac {d}{dt}(g.exptX.g^{-1}.gx)\mid
_{t=0}=dg_x(\overline{X_x}).\] Hence, for $X,Y\in \cal {G}$,$x\in
M$ and $g\in G$, we get (using the fact that $G$ acts on $M$ by
isometries) :
\begin{eqnarray*}
             g.\Phi_{x}(X,Y)&=h_{x}(\overline {Ad_{g^{-1}}X}_{x},
             \overline {Ad_{g^-1}Y}_{x})&=h_{x}(dg^{-1}_{gx} \overline
{X}_{gx}, dg^{-1}_{gx} \overline {Y}_{gx}) \\
                         &=h_{gx}(\overline {X}_{gx},
\overline{Y}_{gx})& \\
                         &=\Phi_{gx}(X, Y).
             \end{eqnarray*}
Observe that if $G$ acts non-properly on $M$, then so does it on
$\Phi(M)$.

 \subsection{Root  decomposition}

Let $\cal A$ be  a Cartan subalgebra, that is,  a maximal abelian ${\Bbb R}$-split  subalgebra of
$\cal G$ and $A$ the associated Cartan group.  Let $\Phi=\Phi\left({\cal A},\ {\cal G}\right)$ be  the  root system  of
$\left(\cal A,\ \cal G\right)$ and
 \[{\cal G}={\cal G}_{0}\oplus \bigoplus_{\alpha \in
\Phi } {\cal G}_{\alpha}\]  the root space decomposition where
\[{\cal G}_{\alpha}=\{X\in {\cal G}: adA.X=\alpha (A).X ,
\forall A \in {\cal A}\}\]
 \[ {\cal G}_{0}=\{X\in {\cal G}: adA.X=0 , \forall A \in
{\cal A}\}.\]
 Then $A$ acts on $\cal {G}$ by diagonal matrices, since
 \[Ad_{g^{-1}}=Ad_{expH}=e^{ad{H}}=diag(e^{\alpha(H)})_{\alpha \in
 \phi \cup \{0\}}\]
 where $g^{-1}=exp(H)$,  $H \in \cal{A}$. It follows that $A$ acts
 by diagonal matrices, and $S^{2}(\cal G)$ admits the following
 decomposition \[S^{2}({\cal G})=\bigoplus_{\lambda \in \Phi \cup
 \{0\}+ \Phi \cup \{0\}} {V_{\lambda}}.\]
 where $V_{\lambda}$ is the set of symmetric bilinear forms $q$ on
 $\cal G$ which satisfy:
\[q(exp(H).X_{1},exp(H).X_{2})= e^{\lambda (H)}.q(X_{1},X_{2}),\]
for all $H\in \cal{A}$ and all $X_{1},X_{2} \in
 \cal G $.
 Keeping in mind that for $X_{1}\in \cal{G_\alpha}$ and $X_{2}\in \cal{G_\beta}$
 we have :  \[q(exp(H).X_{1},exp(H).X_{2})= e^{(\alpha
 +\beta)(H)}.q(X_{1},X_{2}).\] It follows that the forms $q$ in
 $V_{\lambda}$ satisfy :
 \[\alpha +\beta\neq\lambda\Rightarrow
 \cal{G_{\alpha}}\bot \cal{G_{\beta}}.\]

\subsection{Properness of abelian actions}
The following is a  criterion for the non-properness of
linear actions of abelian  Lie groups
 \begin{lemma}
 Let $\{\lambda_1,\cdots \lambda_n\}$ be a generating system in $\Bbb{R}^d$. Let $\Bbb{R}^d$ act faithfully on
$\Bbb{R}^n$ by diagonal matrices as follows. For  $t\in \Bbb{R}^d$, set
 $M(t)=diag(e ^{\langle\lambda_{i},t\rangle})_{1\leq i \leq n}$, where
 $\langle .,.\rangle$ is the usual inner product in $\Bbb{R}^d$. Assume $V$ is an invariant
 (topological) subspace of $\Bbb{R}^n$ on which the action is nonproper. Then there exists
a nonzero vector $t_0\in \Bbb{R}^d$ and an element $x\in
\Bbb{R}^n$ such that $x_i=0$ if $\lambda_i (t_0)< 0$ or an element
$y\in \Bbb{R}^n$ such that $y_i=0$ if $\lambda_i (t_0)> 0$.
\end{lemma}
\begin{proof}
 Since the action  on $V$ is nonproper,  there
exists a sequence $(t_{p})$ with $ t_{p}\rightarrow +\infty$ in
$\Bbb{R}^d$ and a sequence $(x_{p})$ in $V$  such that
$x_{p}\rightarrow x$ in $V$ and $y_{p}=t_{p}.x_{p}\rightarrow y$
in $V$. Consider the sequence $\frac{t_p}{\| t_p\|}$. Up to taking
a subsequence, we may assume it has a limit $t_0$. Since the
action is faithful, and $t_0\neq 0$, there exists $i\in
{1,\cdots,d}$ such that $\lambda_i(t_0)\neq 0$. Note that
$\lambda_i(t_p)\rightarrow + \infty$ if $\lambda_i(t_0)> 0$ and
$\lambda_i(t_p)\rightarrow -\infty$ if $\lambda_i(t_0)< 0$. Hence
$x^i= 0$ if $\lambda_i(t_0)> 0$ and $y^i= 0$ if $\lambda_i(t_0)<
0$.
 \end{proof}
\subsection{End of the proof}
 As we mentioned above, $G$ acts nonproperly on $\Phi(M)$. Let $G=KAK$ be the cartan decomposition
of $G$. Since $G$ has finite center $K$ is compact. So $A$ acts
also nonproperly on $\Phi(M)$. From this, it follows that there
exists $t\neq 0$, $q\in \Phi(M)$ and $\lambda_{0} \in \Phi$ such
that $\lambda_{0}(t)< 0$ and $q_{\lambda}=0$ for all $\lambda \in
\Phi$ with $\lambda(t)< 0$. Put $q=\Phi_{x}$. Then
$\bigoplus_{\alpha(t)<0}\cal{G}_\alpha$ is isotropic with respect
to $\Phi_x$. Hence the image of
$\bigoplus_{\alpha(t)<0}\cal{G}_\alpha$ by the map
$X\mapsto\overline {X}_{x}$ is an isotropic subspace of $T_{x}M$,
so its dimension is less or equal to $1$.  However:

\begin{fact}
 For any $t \in {\cal A}$, the dimension of
$\bigoplus_{\alpha(t)<0}\cal{G}_\alpha$ is at least 2  (where $G$
is assumed to have no local factor locally  isomorphic to
$SL(2, {\Bbb R})$).
\end{fact}
\begin{proof}
 This dimension can not be 0, since $G$ is semi-simple. If it equals
 1, then,  the subalgebra $\bigoplus_{\alpha(t) \geq 0}\cal{G}_\alpha$ has codimension 1 in ${\cal G}$. This contradicts the non existence of $SL(2, {\Bbb R})$ factor condition
 (only simple groups locally isomorphic  to $SL(2, {\Bbb R})$ act on $1$-manifolds).
\end{proof}

We infer from this the existence of
 a nonzero element $X\in  \bigoplus_{\alpha(t)<0}\cal{G}_\alpha$
 such
that $\overline {X}_{x}=0$, which yields : $exp(tX)\in stab(x),
\forall t \in \Bbb R$.  But elements of  $\bigoplus_{\alpha(t)<0}\cal{G}_\alpha$
are nilpotent, and thus generate non-compact groups. This finishes the proof of Theorem
\ref{action.homogeneous}.
 \hfill $\diamondsuit$


\medskip
\noindent
A. Arouche, M. Deffaf, \\
 Facult\'e des math\'ematiques, USTHB \\
BP 32 El'Alia, Bab Ezzouar, Alger, Algeria. \\
 
\noindent
A. Zeghib, \\
CNRS, UMPA, \'Ecole Normale Sup\'erieure de Lyon \\
46, all\'ee d'Italie,
 69364 Lyon cedex 07,  France \\ \\
arouche@math.usthb.dz \\
deffaf1@yahoo.fr \\
Zeghib@umpa.ens-lyon.fr, www.umpa.ens-lyon.fr/\~{}zeghib/
\end{document}